\documentclass[12pt]{article}
\usepackage{amsfonts}

\begin {document}
\topmargin= -.2in \baselineskip=20pt

\title {Mirror Congruence For Rational Points On Calabi-Yau Varieties}

\author {Lei Fu\\
{\small Institute of Mathematics, Nankai University, Tianjin, P.
R. China}\\
{\small leifu@nankai.edu.cn}\\{}\\
Daqing Wan\\
{\small Institute of Mathematics, Chinese Academy of Sciences,
Beijing, P.R. China}
\\
{\small Department of Mathematics, University of California,
Irvine, CA
92697}\\
{\small dwan@math.uci.edu}}
\date{}
\maketitle

\bigskip
\bigskip
\centerline {\bf 0. Introduction}

\bigskip
\bigskip
One of the basic problems in arithmetic mirror symmetry is to
compare the number of rational points on a mirror pair of
Calabi-Yau varieties. At present, no general algebraic geometric
definition is known for a mirror pair. But an important class of
mirror pairs comes from certain quotient construction. In this
paper, we study the congruence relation for the number of rational
points on a quotient mirror pair of varieties over finite fields.
Our main result is the following theorem:

\bigskip
\noindent {\bf Theorem 0.1.} Let $X_0$ be a smooth projective
variety over the finite field ${\bf F}_q$ with $q$ elements of
characteristic $p$. Suppose $X_0$ has a smooth projective lifting
$X$ over the Witt ring $W=W({\bf F}_q)$ such that the $W$-modules
$H^r(X,\Omega_{X/W}^s)$ are free. Let $G$ be a finite group of
$W$-automorphisms acting on the right of $X$. Suppose $G$ acts
trivially on $H^i(X, {\cal O}_{X})$ for all $i$. Then for any
natural number $k$, we have the congruence
$$\# X_0({\bf F}_{q^k})\equiv \# (X_0/G) ({\bf F}_{q^k})~({\rm mod}~q^k),$$
where $\# X_0({\bf F}_{q^k})$ (resp. $\#(X_0/G)({\bf F}_{q^k})$)
denotes the number of elements of the sets of ${\bf
F}_{q^k}$-rational points of $X_0$ (resp. $X_0/G$).

\bigskip
The main application of the above theorem is to Calabi-Yau
varieties. This gives the following theorem announced in [W],
which was the main motivation of the present paper.

\bigskip
\noindent {\bf Theorem 0.2.} Let $X_0$ be a geometrically
connected smooth projective Calabi-Yau variety of dimension $n$
over the finite field ${\bf F}_q$ with $q$ elements of
characteristic $p$. Suppose $X_0$ has a smooth projective lifting
$X$ over the Witt ring $W=W({\bf F}_q)$ such that the $W$-modules
$H^r(X,\Omega_{X/W}^s)$ are free. Let $G$ be a finite group of
$W$-automorphisms acting on the right of $X$. Suppose $G$ fixes a
non-zero $n$-form on $X$. Then for any natural number $k$, we have
the congruence
$$\# X_0({\bf F}_{q^k})\equiv \# (X_0/G) ({\bf F}_{q^k})~({\rm mod}~q^k).$$

\bigskip
\noindent {\bf Proof.} If $X$ is a Calabi-Yau scheme over $W$ of
dimension $n$, then $H^i(X,{\cal O}_X)=0$ for $i\not = 0, n$ and
$G$ acts trivially on them. If the generic fiber of $X$ is
geometrically connected, then $G$ acts trivially on $H^0(X,{\cal
O}_X)$. By Serre duality, $H^n(X,{\cal O}_X)$ is dual to $H^0(X,
\Omega_{X/W}^n)$. Since $X$ is Calabi-Yau, $\Omega_{X/W}^n$ is a
trivial invertible sheaf. In order for $G$ to act trivially on
$H^n(X,{\cal O}_{X/W})$, it suffices for $G$ to fix a nonzero
$n$-form. Theorem 0.2 thus follows from Theorem 0.1.

\bigskip
In particular, we have the following corollary:

\bigskip
\noindent {\bf Corollary 0.3.} Let $X_0$ be the smooth
$(n-1)$-dimensional hypersurface
$$x_0^{n+1}+\cdots +x_n^{n+1}+\lambda x_0\cdots x_n=0$$
in ${\bf P}^n_{{\bf F}_q}$, where $\lambda\in {\bf F}_q$. Let
$$G=\{(\zeta_0,\ldots, \zeta_n)|\zeta_i\in {\bf F}_q,
\zeta_i^{n+1}=1, \prod_{i=0}^n\zeta_i=1\}.$$  Consider the action
$G\times X_0\to X_0$ defined by
$$(\zeta_0,\ldots,\zeta_n)\times [x_0:\ldots: x_n]\mapsto
[\zeta_0x_0: \ldots:\zeta_nx_n].$$ We have $\#X_0({\bf
F}_{q^k})\equiv \#(X_0/G)({\bf F}_{q^k})~({\rm mod}~q^k)$ for any
natural number $k$.

\bigskip
It is well known that the above hypersurface is Calabi-Yau. A
$G$-equivariant nonzero $(n-1)$-form is $\frac {(-1)^i
dx_0\wedge\cdots\wedge \widehat {dx_i} \wedge\cdots \wedge
dx_n}{1+\sum_{j\not =i}x_j^{n+1}-\lambda \prod_{j\not =i} x_j}$ on
the affine space  $x_i =1$ of ${\bf P}^n$.

\bigskip
It is known that for the above hypersurface $X_0$, $X_0/G$ is a
strong singular mirror of $X_0$ if $(n+1)|(q-1)$. It is
conjectured in [W] that for a strong mirror pair of Calabi-Yau
varieties $\{X_0, X_0^{\prime}\}$ over the finite field ${\bf
F}_q$, we have $\#X_0({\bf F}_{q^k})\equiv \# X_0^\prime ({\bf
F}_{q^k})~({\rm mod}~q^k)$ for any integer $k$. See [W] for a
fuller discussion on this and other arithmetic mirror conjectures.
In the situation of Theorem 0.2, if $X/G$ is a singular mirror of
$X$ and if $Y$ is a smooth crepant resolution of $X/G$, then the
pair $(X,Y)$ forms a strong mirror pair of smooth projective
Calabi-Yau varieties. The congruence mirror conjecture in this
case then reduces to showing the congruence
$$\#(X/G)({\bf F}_{q^k})\equiv \# Y({\bf F}_{q^k})~({\rm mod}~q^k).$$

\bigskip
Another application of the theorem is to geometrically connected
varieties with the property $H^i(X,{\cal O}_X)=0$ for all $i\not =
0$. Again in this case, $G$ acts trivially on $H^i(X,{\cal O}_X)$
for all $i$. Let $\overline K$ be the algebraic closure of the
fraction field of $W=W({\bf F}_q)$. By [E], if the $l$-adic
cohomology group $H^i(X\otimes _W \overline K,{\bf Q}_l)$
satisfies the coniveau 1 condition for each $i\not =0$, that is,
if any cohomology class in $H^i(X\otimes _W \overline K,{\bf
Q}_l)$ vanishes in $H^i(U, {\bf Q}_l)$ when restricted to some
nonempty open $U\subset X\otimes _W \overline K$, then we have
$H^i(X,{\cal O}_X)=0$ for all $i\not = 0$. The converse is true if
we assume the generalized Hodge conjecture. It turns out that in
this case, we can prove a theorem stronger than Theorem 0.1. We
don't need to assume $X_0$ can be lifted to $W$.

\bigskip
\noindent {\bf Theorem 0.4.} Let $X_0$ be a smooth geometrically
connected projective variety over the finite field ${\bf F}_q$.
Suppose $H^i(X_0, {\cal O}_{X_0})=0$ for all $i\not =0$. Then for
any natural number $k$, we have
$$\# X_0({\bf F}_{q^k})
\equiv 1 ~({\rm mod}~q^k).$$ Let $G$ be a finite group of ${\bf
F}_q$-automorphisms acting on the right of $X_0$. We have
$$\# (X_0/G) ({\bf F}_{q^k}) \equiv \# X_0({\bf F}_{q^k}) \equiv 1
~({\rm mod}~q^k).$$

\bigskip
\noindent {\bf Acknowledgements.} Our proof is based on
crystalline cohomology and the Mazur-Ogus theorem. H. Esnault
informed us that the results of the present paper can also be
derived using de Rham-Witt cohomology and rigid cohomology. The
research of Lei Fu is supported by the Qiushi Science \&
Technologies Foundation, by the Fok Ying Tung Education
Foundation, by the Transcentury Training Program Foundation, by
the Project 973, and by the SRFDP. The research of Daqing Wan is
partially supported by NSF. Part of this work is done while Lei Fu
is visiting the University of California at Irvine. He would like
to thank the Mathematics Department for its hospitality.

\bigskip
\bigskip
\centerline {\bf 1. Proof of the Theorems}

\bigskip
\bigskip
First we introduce some notations. For any smooth proper scheme
$X_0$ over ${\bf F}_q$, let $H^i(X_0/W)$ be the crystalline
cohomology group of $X_0$. It is a finitely generated module over
the Witt ring $W=W({\bf F}_q)$. Denote by $F:X_0\to X_0$ the
Frobenius correspondence, that is, it is the identity map on the
underlying topological space of $X_0$, and it maps a section of
${\cal O}_{X_0}$ to its $q$-th power.

Let $\kappa$ be a field and let $Z$ be a scheme over $\kappa$.
Denote by $|Z|$ the set of Zariski closed points in $Z$. For any
$z\in |Z|$, define ${\rm deg}(z)=[k(z):\kappa]$, where $k(z)$ is
the residue field at $z$. Let $f:Z\to Z$ be a
$\kappa$-endomorphism with isolated fixed points. Set
$$Z^f=\{z\in |Z| | f(z)=z \hbox { and } f \hbox { induces identity on }k(z)\},$$
and define
$$\Lambda(f)=\sum_{z\in Z^f} {\rm deg}(z).$$
Let $\kappa^\prime$ be a field extending $\kappa$ and let
$f^\prime: Z\otimes_\kappa \kappa^\prime\to Z\otimes_\kappa
\kappa^\prime$ be the base change of $f$. Then we have
$\Lambda(f)=\Lambda(f^\prime)$.

\bigskip
\noindent {\bf Lemma 1.1.} Let $X_0$ be a smooth projective
variety over the finite field ${\bf F}_q$, let $g:X_0\to X_0$ be
an ${\bf F}_q$-automorphism of finite order, and let $K={\rm Frac}
W$ be the fraction field of $W=W({\bf F}_q)$. Then ${\rm
Tr}(F^k,H^i(X_0/W)\otimes_W K)$ and ${\rm
Tr}(gF^k,H^i(X_0/W)\otimes_W K)$ are algebraic integers for any
positive integer $k$ and any $i$, and
\begin{eqnarray*}
\Lambda(F^k)&=&\sum_{i=0}^{2{\rm dim}X_0} (-1)^i{\rm Tr}(F^k,
H^i(X_0/W)\otimes_W K), \\
\Lambda(gF^k)&=&\sum_{i=0}^{2{\rm dim}X_0} (-1)^i{\rm Tr}(gF^k,
H^i(X_0/W)\otimes_W K).
\end{eqnarray*}

\bigskip
\noindent {\bf Proof.} Let $l$ be a prime number distinct from
$p$. By Deligne's theorem ([D] 3.3.9), ${\rm
Tr}(F^k,H^i(X_0\otimes _{{\bf F}_q}\overline {{\bf F}}_q,
\overline {\bf Q}_l))$ are algebraic integers. By the comparison
theorem of Katz-Messing ([KM]), we have $${\rm
Tr}(F^k,H^i(X_0/W)\otimes_W K)={\rm Tr}(F^k,H^i(X_0\otimes _{{\bf
F}_q}\overline {{\bf F}}_q, \overline {\bf Q}_l)).$$ So ${\rm
Tr}(F^k,H^i(X_0/W)\otimes_W K)$ are algebraic integers.  The
formula for $\Lambda(F^k)$ follows from the Lefschetz fixed point
formula in crystalline cohomology theory ([B] Th\'eor\`eme VII
3.1.9).

We will reduce the statements about $gF^k$ to the corresponding
statements for $F^k$. Suppose $g:X_0\to X_0$ has finite order $m$.
Let $X_1=X_0\times_{{\rm Spec}{\bf F}_q}{\rm Spec}{\bf F}_{q^m}$,
and let $\varphi\in {\rm Gal}({\bf F}_{q^m}/{\bf F}_q)$ be the
Frobenius substitution. For any $\sigma\in {\rm Gal}({\bf
F}_{q^m}/{\bf F}_q)$, we have $\sigma=\varphi^k$ for some integer
$k$ uniquely determined modulo $m$. Define
$$f_\sigma:X_1\to
X_1$$ to be the isomorphism of schemes
$$f_\sigma= ({\rm id}_{X_0}\times
\sigma^\ast)\circ (g^{-k}\times {\rm id}_{{\rm Spec}{\bf
F}_{q^m}}):X_0\times_{{\rm Spec}{\bf F}_q}{\rm Spec}{\bf
F}_{q^m}\to X_0\times_{{\rm Spec}{\bf F}_q}{\rm Spec}{\bf
F}_{q^m}.$$ Note that $f_\sigma$ is independent of the choice of
$k$ since $g$ has order $m$. Since $g^{-k}\times {\rm id}_{{\rm
Spec}{\bf F}_{q^m}}$ is an ${\bf F}_{q^m}$-morphism of $X_1$, the
following diagram commutes:
$$\begin{array}{ccc}
X_1&\stackrel {f_\sigma}\to & X_1\\
\downarrow&&\downarrow \\
{\rm Spec}{\bf F}_{q^m}&\stackrel{\sigma^\ast}\to& {\rm Spec}{\bf
F}_{q^m}.
\end{array}
$$ Moreover we have $$f_\tau f_\sigma=f_{\sigma\tau}$$
for any $\sigma,\tau\in {\rm Gal}({\bf F}_{q^m}/{\bf F}_q)$. By
the theory of galois descent, ([S] Chapter V, No. 20, or
Corollarie 7.7 in [SGA 1] Expos\'e VIII), there exists a scheme
$X_0^\prime$ over ${\rm Spec}{\bf F}_q$ such that we have an ${\bf
F}_{q^m}$-isomorphism $$X_1\cong X_0^\prime\times_{{\rm Spec}{\bf
F}_q}{\rm Spec}{\bf F}_{q^m}$$ and the following diagrams commute:
$$\begin{array}{rcl}
X_1&\stackrel {f_\sigma}\to& X_1\\
\cong\downarrow&&\downarrow\cong\\
X_0^\prime\times_{{\rm Spec}{\bf F}_q}{\rm Spec}{\bf
F}_{q^m}&\stackrel{{\rm id}_{{X_0}^\prime}\times \sigma^\ast} \to&
X_0^\prime\times_{{\rm Spec}{\bf F}_q}{\rm Spec}{\bf F}_{q^m}.
\end{array}$$
For any scheme $Z$ of characteristic $p$, let $F_Z:Z\to Z$ be the
Frobenius correspondence, that is, $F_Z$ is identity on the
underlying topological space and the morphism of sheaves
$F_Z^\sharp: {\cal O}_Z\to F_{Z \ast} {\cal O}_Z$ maps each
section to its $q$-th power. On $X_1=X_0\times_{{\rm Spec}{\bf
F}_q} {\rm Spec}{\bf F}_{q^m}$, we have
\begin{eqnarray*}
F_{X_1}&=&({\rm id}_{X_0}\times \varphi^\ast)\circ (F_{X_0}\times
{\rm id}_{{\rm Spec}{\bf F}_{q^m}})=f_{\varphi}\circ (g\times {\rm
id}_{{\rm Spec}{\bf F}_{q^m}})\circ (F_{X_0}\times {\rm id}_{{\rm
Spec}{\bf
F}_{q^m}})\\
&=&f_\varphi\circ (gF_{X_0}\times {\rm id}_{{\rm Spec}{\bf
F}_{q^m}}).
\end{eqnarray*}
Through the isomorphism $X_1\cong X_0^\prime\times_{{\rm Spec}
{\bf F}_q}{\rm Spec}{\bf F}_{q^m}$, $F_{X_1}$ is identified with
$({\rm id}_{X_0^\prime}\times \varphi^\ast)\circ (F_{X_0^\prime
}\times {\rm id}_{{\rm Spec}{\bf F}_{q^m}})$. Moreover, the
commutative diagram above shows that $f_\varphi$ is identified
with ${\rm id}_{X_0^\prime}\times \varphi^\ast$. So the morphism
$gF_{X_0}\times {\rm id}_{{\rm Spec}{\bf F}_{q^m}}$ on
$X_0\times_{{\bf F}_q}{\bf F}_{q^m}$ is identified with the
morphism $F_{X_0^\prime}\times {\rm id}_{{\rm Spec}{\bf F}_{q^m}}$
on $X_0^\prime \times _{{\rm Spec}{\bf F}_q}{\bf F}_{q^m}$. So we
have
\begin{eqnarray*}
&&{\rm Tr}\Biggl(gF_{X_0}\times {\rm id}_{{\bf F}_{q^m}},
H^i\biggl(X_0\times _{{\bf F}_q}{\bf F}_{q^m}/W({\bf
F}_{q^m})\biggr)\otimes_{W({\bf
F}_{q^m})} {\rm Frac} (W({\bf F}_{q^m}))\Biggr)\\
&=&{\rm Tr}\Biggl(F_{X_0^\prime}\times {\rm id}_{{\bf F}_{q^m}},
H^i\biggl(X_0^\prime \times _{{\bf F}_q}{\bf F}_{q^m}/W({\bf
F}_{q^m})\biggr)\otimes_{W({\bf F}_{q^m})} {\rm Frac} (W({\bf
F}_{q^m}))\Biggr).
\end{eqnarray*}
By the base change theorem in crystalline cohomology theory ([B]
Corollaire V 3.5.7), we have
\begin{eqnarray*}
&&{\rm Tr}(gF_{X_0}, H^i(X_0/W)\otimes_W K)\\&=&{\rm
Tr}\Biggl(gF_{X_0}\times {\rm id}_{{\bf F}_{q^m}},
H^i\biggl(X_0\times _{{\bf F}_q}{\bf F}_{q^m}/W({\bf
F}_{q^m})\biggr)\otimes_{W({\bf F}_{q^m})} {\rm
Frac} (W({\bf F}_{q^m}))\Biggr),\\
&&{\rm Tr}(F_{X_0^\prime}, H^i(X_0^\prime /W)\otimes_W K)\\&=&{\rm
Tr}\Biggl(F_{X_0^\prime}\times {\rm id}_{{\bf F}_{q^m}},
H^i\biggl(X_0^\prime \times _{{\bf F}_q}{\bf F}_{q^m}/W({\bf
F}_{q^m})\biggr)\otimes_{W({\bf F}_{q^m})} {\rm Frac} (W({\bf
F}_{q^m}))\Biggr).
\end{eqnarray*}
So we have $${\rm Tr}(gF_{X_0}, H^i(X_0/W)\otimes_W K)={\rm
Tr}(F_{X_0^\prime}, H^i(X_0^\prime /W)\otimes_W K).$$ In
particular, ${\rm Tr}(gF_{X_0}, H^i(X_0/W)\otimes_W K)$ are
algebraic integers for all $i$. Moreover, we have
\begin{eqnarray*}
\Lambda(gF_{X_0})&=&\Lambda(gF_{X_0}\times {\rm id}_{{\rm
Spec}{\bf F}_{q^m}})\\
&=&\Lambda(F_{{X_0}^\prime}\times {\rm id}_{{\rm
Spec}{\bf F}_{q^m}})\\
&=& \Lambda(F_{X_0}^\prime) \\
&=&\sum_{i=0}^{2{\rm dim}X_0} (-1)^i {\rm Tr}(F_{X_0^\prime},
H^i(X_0^\prime /W)\otimes_W K)\\
&=& \sum_{i=0}^{2{\rm dim}X_0} (-1)^i {\rm Tr}(gF_{X_0},
H^i(X_0/W)\otimes_W K).
\end{eqnarray*}
This proves the statements for $gF$. To prove the statements for
$gF^k$, we use the base change from ${\bf F}_q$ to ${\bf
F}_{q^k}$.

\bigskip
\noindent {\bf Lemma 1.2.} Under the condition of Theorem 0.1, we
have
$${\rm Tr}(gF^k, H^i(X_0/W)\otimes_W K)\equiv
{\rm Tr}(F^k, H^i(X_0/W)\otimes_W K)\;({\rm mod}\; q^k)$$ for all
$i$.

\bigskip
\noindent {\bf Proof.} Let $H^i=H^i(X_0/W)$. Recall that $H^i$ can
be identified with the de Rham cohomology of the lifting $X$ of
$X_0$ to $W=W({\bf F}_q)$. (Confer [B] Th\'eor\`eme V 2.3.2). On
$H^i$, we have the Hodge filtration
$$H^i=F^0H^i\supset F^1H^i\supset\cdots$$ and this
filtration is $G$ stable. By a result of Mazur (the property (8.2)
on page 65 of [M]), we have
$$F(F^1H^i)\subset q H^i.$$
We have
$$H^i/F^1H^i=F^0H^i/F^1H^i\cong H^i(X, {\cal O}_X).$$
Choose a basis $\{e_1,\ldots, e_s\}$ of $F^1H^i$ and extend it to
a basis $\{e_1,\ldots, e_s,e_{s+1},\ldots,e_{s+t}\}$ of $H^i$.
Since $F^k(F^1H^i)\subset q^kH^i$, the matrix of $F^k$ on $H^i$
with respect to the above basis is of the form
$$\left(
\begin{array}{cc}
q^kA&q^kB\\
C&D
\end{array}
\right),$$ where $A$ is an $s\times s$ matrix, $B$ is an $s\times
t$ matrix, $C$ is a $t\times s$ matrix, and $D$ is a $t\times t$
matrix. Since $G$ acts trivially on $H^i/F^1H^i\cong H^i(X,{\cal
O}_X)$ and $G$ preserves the Hodge filtration, the matrix of $g\in
G$ on $H^i$ with respect to the above basis is of the form
$$\left(
\begin{array}{cc}
P&O\\
Q&I
\end{array}
\right),$$ where $P$ is an $s\times s$ matrix, $O$ is the $s\times
t$ zero matrix, $Q$ is a $t\times s$ matrix, and $I$ is the
$t\times t$ identity matrix. So the matrix of $gF^k$ is
$$\left(
\begin{array}{cc}
q^kA&q^kB\\
C&D
\end{array}
\right)\left(
\begin{array}{cc}
P&O\\
Q&I
\end{array}
\right)=\left(
\begin{array}{cc}
q^kAP+q^kBQ&q^kB\\
CP+DQ&D
\end{array}
\right).$$ We have $${\rm Tr}(gF^k, H^i)={\rm
Tr}(q^kAP+q^kBQ)+{\rm Tr}(D).$$ On the other hand, we have $${\rm
Tr}(F^k,H^i)={\rm Tr}(q^kA)+{\rm Tr}(D).$$ So we have
$${\rm Tr}(gF^k, H^i)\equiv {\rm Tr}(F^k, H^i)\; ({\rm mod} \;q^k).$$
This finishes the proof of Lemma 1.2.

\bigskip
\noindent {\bf Lemma 1.3.} Let $X_0$ be a quasi-projective scheme
over ${\bf F}_q$, let $G$ be a finite group acting on the right of
$X_0$. Then for any natural number $k$, we have
$$\#(X_0/G)({\bf F}_{q^k})=\frac{1}{\#G} \sum_{g\in G}
\Lambda(gF^k).$$

\bigskip
\noindent {\bf Proof.} This result is well known. We include a
proof here for completeness. Let $Y_0=X_0/G$, and let $|X_0|$
(resp. $|Y_0|$) be the set of Zariski closed point in $X_0$ (resp.
$Y_0$). For any $x\in |X_0|$, define the decomposition subgroup at
$x$ by
$$G_d(x)=\{g\in G| gx=x\}$$
and the inertia subgroup at $x$ by
$$G_i(x)=\{g\in G_d(x)| g \hbox{ {\rm induces identity on the residue
field} } k(x) \hbox{ {\rm at} } x\}.$$ Let $y$ be the image of $x$
in $Y_0$. By Proposition 1.1 in Expos\'e V of [SGA 1], we have an
isomorphism
$$G_d(x)/G_i(x)\cong {\rm Gal}(k(x)/k(y)),$$ and for any $y\in |Y_0|$,
there are exactly $\frac{\#G}{\#G_d(x)}$ Zariski closed points in
$X_0$ above $y$ and each of these closed points has degree ${\rm
deg}(y)\frac{\#G_d(x)}{\#G_i(x)}$. We have
\begin{eqnarray*}
\#Y_0({\bf F}_{q^k})&=& \sum_{y\in |Y_0|,\; {\rm deg}(y)|k} {\rm
deg}(y)\\
&=&\frac{1}{\#G}\sum_{y\in |Y_0|, \;{\rm deg}(y)|k}
\frac{\#G}{\#G_d(x)} \frac{\#G_d(x)}{\#G_i(x)}\#G_i(x){\rm
deg}(y)\\
&=&\frac{1}{\#G} \sum_{y\in |Y_0|,\;{\rm deg}(y)|k}\;\sum_{x\in
|X_0|,\;x\mapsto y} {\rm deg}(x)\#G_i(x).
\end{eqnarray*}
Let $y\in |Y_0|$ be a Zariski closed point with ${\rm deg}(y)|k$,
let $x\in |X_0|$ be a point above $y$, and let $\phi_y\in {\rm
Gal}(k(x)/k(y))$ be the Frobenius substitution. Suppose $g\in
G_d(x)$ and $g^{-1}\mapsto\phi_y^{\frac {k}{{\rm deg}(y)}}$ under
the canonical homomorphism $G_d(x)\to {\rm Gal}(k(x)/k(y))$. Then
$gF^k(x)=x$ and $gF^k$ induces identity on $k(x)$. Conversely, if
$x$ is a Zariski closed point in $X_0$ such that $gF^k(x)=x$ and
$gF^k$ induces identity on $k(x)$, then $g\in G_d(x)$, ${\rm
deg}(y)|k$, and $g^{-1}\mapsto\phi_y^{\frac {k}{{\rm deg}(y)}}$,
where $y$ is the image of $x$ in $Y_0$. On the other hand, there
are exactly $\# G_i(x)$ elements $g$ in $G_d(x)$ such that
$g^{-1}\mapsto\phi_y^{\frac {k}{{\rm deg}(y)}}$. So we finally get
\begin{eqnarray*}
\#Y_0({\bf F}_{q^k})&=&\frac{1}{\#G} \sum_{y\in |Y_0|,\;{\rm
deg}(y)|k}\;\sum_{x\in |X_0|,\;x\mapsto y} {\rm deg}(x)\#G_i(x)\\
&=& \frac{1}{\#G}\sum_{g\in G} \Lambda(gF^k).
\end{eqnarray*}
This proves Lemma 1.3.

\bigskip
Now we are ready to prove Theorem 0.1. By Lemmas 1.3 and 1.1, we
have
\begin{eqnarray*}
\#(X_0/G)({\bf F}_{q^k})&=&\frac{1}{\#G}\sum_{g\in G} \Lambda(gF^k)\\
&=& \frac{1}{\#G}\sum_{g\in G}\sum_{i=0}^{2{\rm dim} X_0}
(-1)^i{\rm Tr}(gF^k, H^i(X_0/W)\otimes_W K).
\end{eqnarray*}
By Lemmas 1.1 and 1.2, ${\rm Tr}(gF^k, H^i(X_0/W)\otimes_W K)$ and
${\rm Tr}(F^k, H^i(X_0/W)\otimes_W K)$ are algebraic integers, and
$${\rm Tr}(gF^k, H^i(X_0/W)\otimes_W K)
\equiv {\rm Tr}(F^k, H^i(X_0/W)\otimes_W K) \; ({\rm mod}\;
q^k).$$  From now on, we work over the integral closure of
$p$-adic integers. Let ${\rm ord}_q(\#G) =c$, a non-negative
rational number. For each $k\geq c$, we have
\begin{eqnarray*}
\#(X_0/G)({\bf F}_{q^k})&=&\frac{1}{\#G}\sum_{g\in
G}\sum_{i=0}^{2{\rm dim} X_0}(-1)^i{\rm Tr}(gF^k,
H^i(X_0/W)\otimes_W K)\\
&\equiv& \frac{1}{\#G}\sum_{g\in G}\sum_{i=0}^{2{\rm dim} X_0}
(-1)^i{\rm Tr}(F^k, H^i(X_0/W)\otimes_W K)
\;({\rm mod}\; q^{k-c})\\
&\equiv & \sum_{i=0}^{2{\rm dim} X_0}(-1)^i{\rm
Tr}(F^k, H^i(X_0/W)\otimes_W K)\;({\rm mod}\; q^{k-c})\\
&\equiv & \#X_0({\bf F}_{q^k})\;({\rm mod}\; q^{k-c}).
\end{eqnarray*}
Let $Z(X_0,T)$ and $Z(X_0,T)$ be the zeta-functions of $X_0$ and
$X_0/G$, respectively. They are rational functions. Recall that we
have
\begin{eqnarray*}
\frac{d}{dT} \log Z(X_0,T)&=&\sum_{k=1}^\infty \# X_0({\bf
F}_{q^k})T^{k-1},\\
\frac{d}{dT} \log Z(X_0/G,T)&=&\sum_{k=1}^\infty \# (X_0/G)({\bf
F}_{q^k})T^{k-1}.
\end{eqnarray*}
Take a factorization
$$\frac{Z(X_0, T)}{Z(X_0/G, T)} =\prod_{i=1}^m
(1-\alpha_iT)^{-n_i},\ \alpha_i \not=0$$ where the $\alpha_i$'s
are distinct and the $n_i$'s are non-zero integers. Taking
logarithmic derivative on both sides, we get
$$\sum_{k=1}^{\infty} (\#X_0({\bf F}_{q^k})-\#(X_0/G)({\bf
F}_{q^k}))T^{k-1} =\sum_{i=1}^m \frac{n_i\alpha_i}{1-\alpha_iT}.$$
Using the congruence
$$\#(X_0/G)({\bf F}_{q^k}) \equiv \#X_0({\bf F}_{q^k})\;({\rm
mod}\; q^{k-c})$$ for all $k\geq c$, one deduces that the above
power series is $p$-adic analytic in the open disk ${\rm
ord}_q(T)>-1$. This implies that each $\alpha_i$ satisfies ${\rm
ord}_q(\alpha_i)\geq 1$, that is, each $\alpha_i$ is divisible by
$q$. We conclude that
$$\#X_0({\bf F}_{q^k})-\#(X_0/G)({\bf
F}_{q^k})= \sum_{i=1}^m n_i\alpha_i^k \equiv 0 \;({\rm mod}\;
q^{k}).$$ This finishes the proof of Theorem 0.1.

\bigskip
Let's prove Theorem 0.4. By Ogus' generalization of Mazur's
theorem ([BO] Theorem 8.39), the Newton polygon of the Frobenius
correspondence $F$ on $H^i(X_0/W)\otimes_W K$ lies on or above the
Hodge polygon of $X_0$. For any $i\not=0$, we have $H^i(X_0,{\cal
O}_{X_0})=0$. So the slope of each line segment on the Newton
polygon is at least $1$, that is, all the eigenvalues of $F^k$ on
$H^i(X_0/W)\otimes_W K$ are divisible by $q^k$ (as $p$-adic
integers). So we have
$${\rm Tr}(F^k, H^i(X_0/W)\otimes_W K)\equiv 0 \; ({\rm
mod}\; q^k)$$ for all $i\not =0$. Since $X_0$ is geometrically
connected, we have
$${\rm Tr}(F^k, H^0(X_0/W)\otimes_W K)=1.$$ So by Lemma 1.1, we have
\begin{eqnarray*}
\# X_0({\bf F}_{q^k})&=& \sum_{i=0}^{2{\rm dim}X_0} (-1)^i{\rm
Tr}(F^k, H^i(X_0/W)\otimes_W K) \\
&\equiv& 1 \;({\rm mod}\; q^k).
\end{eqnarray*}
Now let $G$ be a finite group acting on the right of $X_0$. For
any $g\in G$, since $g$ has finite order, the action of $g$ on
$H^i(X_0/W)\otimes_W K$ is diagonalizable and all its eigenvalues
are roots of unity. Combining with the fact that $F$ commutes with
$g$, we see that all the eigenvalues of $gF^k$ on
$H^i(X_0/W)\otimes_W K$ are also divisible by $q^k$ for any $i\not
=0$. So we have
$${\rm Tr}(gF^k, H^i(X_0/W)\otimes_W K)
\equiv {\rm Tr}(F^k, H^i(X_0/W)\otimes_W K)\equiv 0 \; ({\rm
mod}\; q^k)$$ for all $i\not =0$. Since $X_0$ is geometrically
connected, we have
$${\rm Tr}(gF^k, H^0(X_0/W)\otimes_W K)={\rm Tr}(F^k, H^i(X_0/W)\otimes_W K)
=1.$$ Again let ${\rm ord}_q(\#G) =c$. For each $k\geq c$, by
Lemmas 1.1, 1.3, and the above discussion, we have
\begin{eqnarray*}
\#(X_0/G)({\bf F}_{q^k})&=&\frac{1}{\#G}\sum_{g\in
G}\sum_{i=0}^{2{\rm dim} X_0}(-1)^i{\rm Tr}(gF^k,
H^i(X_0/W)\otimes_W K)\\
&\equiv& \frac{1}{\#G}\sum_{g\in G}\sum_{i=0}^{2{\rm dim} X_0}
(-1)^i{\rm Tr}(F^k, H^i(X_0/W)\otimes_W K)
\;({\rm mod}\; q^{k-c})\\
&\equiv & \sum_{i=0}^{2{\rm dim} X_0}(-1)^i{\rm
Tr}(F^k, H^i(X_0/W)\otimes_W K)\;({\rm mod}\; q^{k-c})\\
&\equiv & \#X_0({\bf F}_{q^k})\;({\rm mod}\; q^{k-c}).
\end{eqnarray*}
As in the proof of Theorem 0.1, this implies that
$$\#(X_0/G)({\bf F}_{q^k})\equiv \#X_0({\bf F}_{q^k})\;({\rm mod}\; q^k
).$$ This finishes the proof of Theorem 0.4.

\bigskip
\bigskip
\noindent {\bf References.}

\bigskip
\bigskip
\noindent [B] P. Berthelot, {\it Cohomologie Cristalline des
Sch\'emas de Carast\'eristique $p>0$}, Lecture Notes in
Mathematics 407, Springer-Verlag 1974.

\bigskip
\noindent [BO] P. Berthelt and A. Ogus, {\it Notes on Crystalline
Cohomology}, Princeton University Press 1978.

\bigskip
\noindent [D] P. Deligne, {\it La Conjecture de Weil II}, Publ.
Math. IHES, 52 (1980), 137-252.

\bigskip
\noindent [E] H. Esnault, {\it Deligne's Integrality Theorem in
Unequal Characteristic and Rational Points over Finite Fields},
Preprint.

\noindent [KM] N. Katz and W. Messing, {\it Some Consequences of
the Riemann Hypothesis for Varieties over Finite Fields}, Invent.
Math. 23 (1974), 73-77.

\bigskip
\noindent [M] B. Mazur, {\it Frobenius and the Hodge Filtration
(estimates)}, Ann. of Math. 98 (1973), 58-95.

\bigskip
\noindent [S] J.-P. Serre, {\it Algebraic Groups and Class
Fields}, translation of the French edition, Springer-Verlag, 1988.

\bigskip
\noindent [SGA 1] A. Grothendieck, {\it Rev\^etements \'Etales et
Groupe Fondamental}, Lecture Notes in Mathematics 224,
Springer-Verlag (1971).

\bigskip
\noindent [W] D. Wan, {\it Mirror Symmetry for Zeta Functions},
arXiv:math.AG/0411464, 21 Nov 2004.

\end{document}